 \magnification=\magstep1 
\newcount\sectno
\newcount\subsectno
\newcount\parno
\newcount\equationno
\newif\ifsubsections
\subsectionsfalse

\def\sectnum{\the\sectno} 
\def\subsectnum{\sectnum\ifsubsections .\the\subsectno\fi} 
\def\parnum{\subsectnum .{\the\parno}}
\def\eqnum{\subsectnum .\the\equationno}

\def\abstract#1\endabstract
{
{\abstractfont
    \baselineskip=9pt
    \leftskip=4pc  \rightskip=4pc
    \bigskip
    \noindent
     ABSTRACT.\ #1
\medskip} 
}
 
\def\thanks[#1]#2\endthanks{\footnote{$^#1$}{\footnotefont\kern-6pt #2}}

\newcount\minutes
\newcount\scratch

\def\timestamp{%
\scratch=\time
\divide\scratch by 60
\edef\hours{\the\scratch}
\multiply\scratch by 60
\minutes=\time
\advance\minutes by -\scratch
\the \month/\the\day$\,$---$\,$\hours:\null
\ifnum\minutes< 10 0\fi
\the\minutes}

\def\today{\ifcase\month\or
January\or February\or March\or April\or May\or June\or
July\or August\or September\or October\or November\or December\fi
\space\number\day,\number\year}

\outer\def\newsection #1.\par{\vskip1.5pc plus.75pc \penalty-250
     \subsectno=0
     \parno=0
     \equationno=0
     \advance\sectno by1
     \leftline{\smalltitlefont \sectnum.\hskip 1pc  #1}
                \nobreak \vskip.75pc\noindent}

\outer\def\newsectiontwoline #1/#2/.\par{\vskip1.5pc plus.75pc \penalty-250
     \subsectno=0
     \parno=0
     \equationno=0
     \advance\sectno by1
     \leftline{\smalltitlefont \sectnum.\hskip 1pc  #1}
     \leftline{\smalltitlefont \hskip 22pt #2}
                \nobreak \vskip.75pc\noindent}

\outer\def\newsubsection #1.\par{\vskip1pc plus.5pc\penalty-250
     \parno=0
     \equationno=0
     \advance\subsectno by1
     \leftline{{\bf \subsectnum}\hskip 1pc  #1.}
                \nobreak \vskip.5pc\noindent}

\def\newpar #1.{\advance \parno by1
     \par
 \medbreak \noindent 
      {\bf \parnum. #1.} \hskip 6pt}

\long\def \newclaim #1. #2\par {\advance \parno by1
    \medbreak \noindent 
     {\bf \parnum \hskip 6 pt #1.\hskip 6pt} {\sl #2} \par \medbreak}

\def\eq $$#1$${\global \advance \equationno by1 $$#1\eqno(\eqnum)$$}

\def\rmarginsert[#1]{\hglue 0pt\vadjust
{\null\vskip -\baselineskip\rightline{\abstractfont\rlap{\hfil\  #1}}}}

\def\lmarginsert[#1]{\hglue 0pt\vadjust
{\null\vskip -\baselineskip\leftline{\abstractfont\llap{#1\ \hfill}}}}

\newif\ifproofmode
\proofmodefalse

\def\refpar[#1]#2.{\advance \parno by1
     \par
 \medbreak \noindent 
      {\bf \parnum \hskip 6 pt #2.\hskip 6pt}%
\expandafter\edef\csname ref#1\endcsname
{\parnum}\ifproofmode\rmarginsert[\string\ref#1]\fi}

\long\def \refclaim[#1]#2. #3\par {\advance \parno by1
    \medbreak \noindent 
{\bf \parnum \hskip 6 pt #2.\hskip 6pt}%
\expandafter\edef\csname ref#1\endcsname
{\parnum}\ifproofmode\rmarginsert[\string\ref#1]\fi
{\sl #3} \par \medbreak}

\def\refer[#1]{%
\expandafter\xdef\csname ref#1\endcsname
{\parnum}\ifproofmode\rmarginsert[\string\ref#1]\fi}

\def\refereq[#1]$$#2$$ {%
\eq$$#2$$%
\expandafter\xdef\csname ref#1\endcsname{(\eqnum)}%
\ifproofmode\rmarginsert[\string\ref#1]\fi
}

\def\refeq{\refereq}

\def \Definition #1\\ {\vskip 1pc \noindent 
      {\bf #1. Definition. \hskip 6pt}\vskip 1pc}

\def\proof{{PROOF.} \enspace}

\def\qedmark{\hbox{\vrule height 4pt width 3pt}}
\def\qedskip{\vrule height 4pt width 0pt depth 1pc}
\def\qed{\penalty 1000\quad\penalty 1000{\qedmark\qedskip}}

\def \a {\alpha}
\def \b {\beta}

\def \d {\delta}
\def \D {\triangle}
\def \e {\epsilon}

\def \g {\gamma}

\def \K {\nabla}
\def \l {\lambda}
\def \L {\Lambda}
\def \n {\,\vert\,}
\def \N {\,\Vert\,}
\def \o {\theta}
\mathchardef\p="011E    

\def \s {\sigma}

\def \w {\omega}
\def \W {\Omega}


\def\Gtwo{{\mathop{{{\mbi G\/}}\kern-.5pt_{{}_2}}}}
\def\Ffour{{\mathop{{{\mbi F\/}}\kern-2.5pt_{{}_4}}}}
\def\Esix{{\mathop{{{\mbi E\/}}\kern-.5pt_{{}_6}}}}
\def\Eseven{{\mathop{{{\mbi E\/}}\kern-.5pt_{{}_7}}}}
\def\Eeight{{\mathop{{{\mbi E\/}}\kern-.5pt_{{}_8}}}}


\def\cg{{\liefont G}}

\def\ck{{\liefont K}}

\def\cm{{\liefont M}}

\def\cn{{\liefont N}}

\def\co{{\liefont O}}

\def\cp{{\liefont P}}

\def\cs{{\liefont S}}

\def\cu{{\liefont U}}

\def\sdp{
\mathop{\hbox{$\raise 1pt\hbox{$\scriptscriptstyle |$}\kern-2.5pt\times$}}
        }

\def\dag{{\raise 1 pt\hbox{{$\scriptscriptstyle \dagger$}}}}

\def\*{\raise 1.5pt \hbox{*}}

\def\tr{\mathop{\rm tr}\nolimits}

\def\ad{\mathop{\rm ad}\nolimits}

\def\id {\mathop{\rm id}\nolimits}

\def \li{\langle}
\def \ri{\rangle}

                                                          
\def\ifundefined#1{\expandafter\ifx\csname#1\endcsname\relax}

\def\cross{\times}
\def\longerrightarrow{-\kern-5pt\longrightarrow}

\def\star{\lower 1pt\hbox{*}}
\def \nulset {
\raise 1pt\hbox{
\hskip -3pt$\not$\kern -0.2pt \raise .7pt\hbox{${\scriptstyle\bigcirc}$}}}
\def \norm|#1|{\Vert#1\Vert}

\def \interior(#1){#1\kern -6pt \raise 7.5pt 
      \hbox{$\scriptstyle \circ$}{}\hskip 2pt}

\def\twist_#1{\kern -.15em\cross\kern -.30 em{}_{{}_{#1}}\kern .07 em}

\font\cmr=cmr10 at 10pt
\font\cmrviii=cmr8
\font\cmrvi=cmr6

\font\cmrXIV=cmr12 at 14 pt
\font\cmrXX=cmr17 at 20 pt
\font\cmrXXIV=cmr17 at 24 pt
\font\cmbxXII=cmbx12
\font\cmbxsl=cmbxsl10
\font\cmbxslviii=cmbxsl10 at 8pt
\font\cmbxslv=cmbxsl10 at 5pt

      \font\tenrm=cmr10 at 10.3 pt
      \font\sevenrm=cmr7 at 7.21 pt
      \font\fiverm=cmr5 at 5.15 pt
      \font\teni=cmmi10 at 10.3 pt
      \font\seveni=cmmi7 at 7.21 pt
      \font\fivei=cmmi5 at 5.15 pt     
      \font\tensy=cmsy10 at 10.3 pt
      \font\sevensy=cmsy7 at 7.21 pt   
      \font\fivesy=cmsy5 at 5.15 pt    
      \font\tenex=cmex10 at 10.3 pt
      \font\tenbf=cmbx10 at 10.3 pt
      \font\sevenbf=cmbx7 at 7.21 pt    
      \font\fivebf=cmbx5 at 5.15 pt  

\def\UseComputerModern   
{
\textfont0=\tenrm \scriptfont0=\sevenrm \scriptscriptfont0=\fiverm
\def\rm{\fam0\tenrm}
\textfont1=\teni \scriptfont1=\seveni \scriptscriptfont1=\fivei
\def\mit{\fam1} \def\oldstyle{\fam1\teni}
\textfont2=\tensy \scriptfont2=\sevensy \scriptscriptfont2=\fivesy
\def\cal{\fam2}
\textfont3=\tenex \scriptfont3=\tenex \scriptscriptfont3=\tenex
\def\it{\fam\itfam\tenit} 
\textfont\itfam=\tenit
\def\sl{\fam\slfam\tensl} 
\textfont\slfam=\tensl
\def\bf{\fam\bffam\tenbf} 
\textfont\bffam=\tenbf \scriptfont\bffam=\sevenbf
\scriptscriptfont\bffam=\fivebf
\def\tt{\fam\ttfam\tentt} 
\textfont\ttfam=\tentt
\def\abstractfont{\cmrviii}
\def\footnotefont{\cmrviii}
\def\tinyfont{\cmrvi}
\def\smalltitlefont{\cmbxXII}  
\def\titlefont{\cmrXIV}
\def\bigtitlefont{\cmrXX}
\def\verybigtitlefont{\cmrXXIV} 
\textfont9=\cmbxsl \scriptfont9=\cmbxslviii \scriptscriptfont9=\cmbxslv
\def\mbi{\fam9}
\cmr
}  

\def\liefont{\cal}

\def \bs {\bigskip}
\def \ms {\medskip}
\def \ss {\smallskip}

\def \ni {\noindent}

\def\enditem{\item{}\par\vskip-\baselineskip\noindent}
\def\ei{\enditem}

   \baselineskip=14 true pt 
   \hsize 37 true pc \hoffset= 22 true pt
   \voffset= 0 true pt
   \vsize  54 true pc

\def\id{\mathop{\rm id}\nolimits}

\def\ni{\noindent}
\def\bs{\bigskip}
\def\ms{\medskip}

\def\ss{\smallskip}

   \baselineskip=14 true pt 
   \hsize 35 true pc \hoffset= 25 true pt
   \vsize  52 true pc
\UseComputerModern
\subsectionsfalse
\font\smalltitlefont=cmbx10 at 11 pt
%


\def\Bibliography
{

\font\TRten=cmr10 at 10 true pt
\font\TIten=cmti10 at 10 true pt
\font\TBten=cmbx10 at 10 true pt

\def\ourindent{\hfil\vskip-\baselineskip}

 \frenchspacing
 \parindent=0pt

 \def \keyfnt{\TRten}
 \def \authornamefnt{\TRten}
 \def \booktitlefnt{\TIten}
 \def \articletitlefnt{\TRten}
 \def \journalnamefnt{\TIten}
 \def \volumefnt{\TBten}
 \def \publishernamefnt{\TRten}
 \def \pagesfnt{\TRten}
 \def \yearfnt{\TRten}
 \def \commentfnt{\TRten}

 \def \bookitem //##1//##2//##3//##4//##5//##6//##7//##8//
      { \goodbreak{\par\hskip-40pt{\keyfnt [##1]}\ourindent{\authornamefnt ##2,}}
              {\booktitlefnt ##3.\/}\thinspace
              {\publishernamefnt ##4,}
              {\yearfnt ##6.}
              {\commentfnt ##8}
       }

\def \b{\bookitem}

\def \articleitem //##1//##2//##3//##4//##5//##6//##7//##8//
      { \goodbreak{\par\hskip-40pt{\keyfnt [##1]}\ourindent{\authornamefnt ##2,}}
              {\articletitlefnt ##3},
              {\journalnamefnt ##4\/}
              {\volumefnt ##5}
              {\hbox{\yearfnt(\hskip -1pt ##6)}},
              {\pagesfnt ##7.}
              {\commentfnt ##8}
       }
\def \a{\articleitem}

\def \preprintitem //##1//##2//##3//##4//##5//##6//##7//##8//
      { \goodbreak{\par\hskip-40pt{\keyfnt [##1]}\ourindent{\authornamefnt ##2,}}
              {\articletitlefnt ##3},
              {\commentfnt ##8}
       }
\def \p{\preprintitem}

   \vskip 1in
   \centerline{References}
   \vskip .5in
}


\def\Reduce#1  
{
\font\stenrm=cmr10 scaled #1
\font\sninerm=cmr9 scaled #1
\font\seightrm=cmr8 scaled #1
\font\ssevenrm=cmr7 scaled #1
\font\ssixrm=cmr6 scaled #1
\font\sfiverm=cmr5 scaled #1

\font\steni=cmmi10 scaled #1
\font\sninei=cmmi9 scaled #1
\font\seighti=cmmi8 scaled #1
\font\sseveni=cmmi7 scaled #1
\font\ssixi=cmmi6 scaled #1
\font\sfivei=cmmi5 scaled #1

\font\stenit=cmti10 scaled #1
\font\snineit=cmti9 scaled #1
\font\seightit=cmti8 scaled #1
\font\ssevenit=cmti7 scaled #1

\font\stensy=cmsy10 scaled #1
\font\sninesy=cmsy9 scaled #1
\font\seightsy=cmsy8 scaled #1
\font\ssevensy=cmsy7 scaled #1
\font\ssixsy=cmsy6 scaled #1
\font\sfivesy=cmsy5 scaled #1

\font\stenbf=cmbx10 scaled #1 
\font\sninebf=cmbx9 scaled #1
\font\seightbf=cmbx8 scaled #1
\font\ssevenbf=cmbx7 scaled #1
\font\ssixbf=cmbx6 scaled #1
\font\sfivebf=cmbx5 scaled #1

\font\stentt=cmtt10  scaled #1
\font\sninett=cmtt9 scaled #1
\font\seighttt=cmtt8 scaled #1

\font\stenex=cmex10 scaled #1

\font\stensl=cmsl10 scaled #1
\font\sninesl=cmsl9 scaled #1
\font\seightsl=cmsl8 scaled #1

\textfont0=\stenrm \scriptfont0=\ssevenrm \scriptscriptfont0=\sfiverm
\def\rm{\fam0\stenrm}
\textfont1=\steni \scriptfont1=\sseveni \scriptscriptfont1=\sfivei
\def\mit{\fam1} \def\oldstyle{\fam1\steni}
\textfont2=\stensy \scriptfont2=\ssevensy \scriptscriptfont2=\sfivesy
\def\cal{\fam2}
\textfont3=\stenex \scriptfont3=\stenex \scriptscriptfont3=\stenex
\def\it{\fam\itfam\tenit} 
\textfont\itfam=\stenit
\def\sl{\fam\slfam\stensl} 
\textfont\slfam=\stensl
\def\bf{\fam\bffam\stenbf} 
\textfont\bffam=\stenbf \scriptfont\bffam=\ssevenbf
\scriptscriptfont\bffam=\sfivebf
\def\tt{\fam\ttfam\stentt} 
\textfont\ttfam=\stentt
\rm
}  


\def\p{\partial}


{\Reduce{1200}
\centerline{\bf Schr\"odinger  flows on Grassmannians}
}
\bs
\centerline { Chuu-Lian Terng\footnote{$^1$}{Research supported
in  part by 
NSF Grant DMS 9626130}
and  Karen Uhlenbeck\footnote{$^2$}{Research supported in part
by  Sid Richardson 
Regents' Chair Funds, University of Texas system}}
\bs\bs\bs

\centerline{Abstract}
\bs
{\Reduce{800}

\ni
The geometric non-linear Schr\"odinger equation (GNLS) on the complex
Grassmannian manifold 
$M$ is the evolution equation on  the space $C(R,M)$ of paths on $M$: 
$$J_\g\g_t=\K_{\g_x}\g_x,$$  where $\K$ is the Levi-Civita connection of the
K\"ahler metric and $J$ is the complex structure.  GNLS is the
Hamiltonian equation for the energy functional on $C(R,M)$
with respect  to the
symplectic form induced from the K\"ahler form on
$M$.   It has a Lax pair that is gauge
equivalent to the Lax pair of the matrix non-linear Schr\"odinger equation (MNLS). 
We construct via gauge transformations an isomorphism from $C(R,M)$ to
the phase space of the MNLS equation so that the GNLS flow corresponds to the
MNLS flow.  The existence of
global solutions to the Cauchy problem for GNLS and the hierarchy of commuting
flows follows from the correspondence.  Direct geometric constructions show the
flows are given by geometric partial differential equations, and the space of
conservation laws has a structure of a non-abelian Poisson group.  We also
construct a hierarchy of symplectic structures for GNLS. Under pullback, the
known order $k$ symplectic structures correspond to the order $k-2$ symplectic
structures that we find.  The shift by two is a surprise, and  is due to the fact
that the group structures depend on gauge choice.  }

\vfil\eject

\newsection Introduction.\par

Harmonic maps from one Riemannian manifold $N$ to another $M$ are critical
points for the energy 
$$E(\g)={1\over 2} \int_N \n d\g\n^2 d\mu_N.$$
The Euler-Lagrange equation is 
$$\D_\g \g= \tr(\K(d\g))=0.$$
When $N$ is $S^2$ or a $2$-dimensional torus and $M$ is a symmetric space,
the moduli space of harmonic maps has been successfully studied by many authors
using techniques from integrable systems ([U], [BFPP], [BG]).  

The heat flow for harmonic maps from
$N$ to
$M$ is the gradient flow for the energy and has the form
$$\g_t= \D_\g \g.$$
If $M$ is a K\"ahler manifold, we have a complex structure $J:TM\to
TM$, $J^2=-$id, and the equation 
$$J_\g\g_t=\K_\g\g$$
has the type of a non-linear Schr\"odinger equation.  Very little is known about this
equation except when dim$(N)=1$. See recent work by Chang, Shatah and
Uhlenbeck for a discussion of the radially symmetric case in dimension $2$
([CSU]). 
  For $N=R^n$, $n=1,2, 3$ and
$M=S^2$, this equation is a simplification of the Landau-Lifshitz equation for a
continuous anisotropic magnet.  For $N=R^1$ and $M=S^2$, the GNLS is often
referred to as the continuous isotropic Heisenberg ferromagnetic model 
([FT]). 

Let $(M,g,J)$ be a K\"ahler manifold with 
metric $g$ and complex structure $J$.   The one-dimensional {\it geometric
non-linear Schr\"odinger equation (GNLS)\/} with target
$M$ is the evolution equation on the space $C(R,M)$ of
smooth paths from $R$ to
$M$: 
\refeq[aa]$$J_\g\g_t= \K_{\g_x} \g_x,$$ where $\K$ is
the Levi-Civita connection of $g$.  The symplectic form
$$\tau_\g(v_1,v_2)=g( J_\g(v_1), v_2)$$  on $M$ induces a natural
symplectic structure $\hat \tau$ on
$C(R,M)$
\refeq[bz]$$\hat \tau_\g(\xi, \eta)=\int_{-\infty}^\infty
g(J_{\g(x)}(\xi(x)),\eta(x)) dx.$$  
Equation \refaa{} is the Hamiltonian equation of the energy functional with
respect to $\hat \tau$.   
The main goal of
this paper is to study the Hamiltonian  theory of the geometric non-linear
Schr\"odinger equation with target the Grassmannians. The theory can be
extended to the more general case of  compact Hermitian symmetric spaces. 

The matrix non-linear Schr\"odinger equation (MNLS) is
\refeq[ad]$$q_t= i(q_{xx} + 2 qq^*q),$$ where $q$ is a map from $R^2$ to
 the space $\cm_{k\times(n-k)}$ of $k\times (n-k)$ complex matrices. This
equation was first studied by Fordy and Kulish in [FK] as a generalization of the NLS
equation.  Note that if  $q$ is a
$1\times 1$ matrix, this is the NLS equation.  When
$k=1$,
$q$ is a vector and this is known as the vector non-linear Schr\"odinger equation. 

The GNLS with target a Grassmanian is known to be gauge equivalent to the matrix
non-linear Schr\"odinger equation. This correspondence
for $S^2=CP^1$ is contained in the classical Hasimoto transform and is described
in detail by Faddeev-Takhtajen ([FT]). We describe this gauge equivalence in general,
and examine the behavior of the hierarchies of symplectic structures under the
gauge equivalence.  We say that a symplectic structure is of order
$k$ if the corresponding Poisson structure is given by an order $k$
integro-differential operator. Each equation has a natural order zero symplectic
structure that arises from a coadjoint orbit.   Under the equivalence, the two order
zero symplectic structures do not correspond to each other.  In fact, we may regard
the gauge change as naturally generating a hierarchy of symplectic structures,
although we would only generate the even order structures in this fashion. 

Hasimoto ([H]) showed that NLS is equivalent to the equation
of da Rios ([dR]):
$$\a_t=\a_x\times \a_{xx},$$
which models the movement of a thin vortex
filament in a viscous liquid. This equation preserves arc length and $\g=\a_x$
satisfies the GNLS with target
$S^2$. Langer and Perline ([LP]) generalize this to MNLS. They give a
geometric realization of the MNLS as an arc length preserving curve evolution in
$R^{k(n-k)}$. 

 The MNLS is a second flow of the $u(n)$-hierarchy, and the Hamiltonian theory
has been described by many authors.  Let $a$ denote the diagonal matrix with
eigenvalues $i/2$ and $-i/2$ and multiplicities
$k$ and $(n-k)$ respectively:
$$a=\pmatrix{{i\over 2}I_k
&0\cr 0& -{i\over 2}I_{n-k}\cr}.$$  Let $u(n)_a$ denote the centralizer of
$a$, and 
$u(n)_a^\perp$ the orthogonal complement of $u(n)_a$ in $u(n)$.  We find
$$\eqalign{u(n)_a&=\{y\in u(n)\n [y,a]=0\}=\left\{\pmatrix{A_1&0\cr
0&A_2\cr}\biggl| \,\, A_1\in u(k), A_2\in u(n-k)\right\},
\cr  u(n)_a^\perp &=\{y\in u(n)\n <y,z>=0 \,\, {\rm for \, all\,}\,\, z\in
u(n)_a\}\cr &=\left\{\pmatrix{0&q\cr -q^*&0\cr}\n q\in \cm_{k\times
(n-k)}\right\}\simeq \cm_{k\times (n-k)}.\cr}$$
Let $\cs(R,u(n)_a^\perp)$ denote the space of smooth maps from $R$ to
$u(n)_a^\perp$ that are in the Schwartz class.
Since $\ad(a)$ is a
skew-adjoint, linear isomorphism on $u(n)_a^\perp$,  the $2$-form
defined by 
\refeq[bh]$$w(v_1,v_2) =\int_{-\infty}^\infty <-\ad(a)^{-1}(v_1),
v_2>dx$$
is symplectic on $\cs(R,u(n)_a^\perp)$.  
The second flow in the AKNS $u(n)$-hierarchy is 
\refeq[bj]$$u_t=(Q_2)_x+ [u, Q_2],$$ where 
\refeq[dv]$$u=\pmatrix{0&q\cr -q^*&0\cr}, \quad Q_2 = \pmatrix{-iqq^*& iq_x\cr
iq^*_x & iq^*q\cr}.$$
The second flow \refbj{} has a Lax pair
 \refeq[ah]$$\left[{\p\over \p x} + a\l + u, \quad {\p\over \p t} + a\l^2 +
u\l + Q_2\right]=0,$$ i.e., $u$ is a solution of \refbj{} if and only if $u$ satisfies
\refah{} for all $\l$.  This is the Hamiltonian equation of the functional
\refeq[bi]$$F_2(u)=\int_{-\infty}^\infty
-{1\over 4}\tr(u_x^2) +{1\over 8} (\tr(u^2))^2dx$$  on
$\cs(R,u(n)_a^\perp)$ with respect to $w$. 
When written in terms of $q$, equation \refbj{}  is the MNLS \refad{},  
$$\eqalign{F_2(q)&= \int_{-\infty}^\infty -{1\over 2}\N q_x\N^2 + {1\over 2} \N
q\N^4 dx,\cr
w(q_1,	q_2)&=\int_{-\infty}^\infty <-iq_1, q_2>dx.\cr}$$  So we shall also
refer equation
\refbj{} as the MNLS.

We list some of the known
properties of the MNLS equation \refbj{}:

\item {(1)}  Beals and Coifman construct the inverse
scattering theory for the first order system $d_x+a\l +u$ to solve
the Cauchy problem globally with Schwartz class initial data ([BC]).   
\item {(2)} We constructed in [TU2] an action of the rational loop group on the
space of solutions of the MNLS equation such that the action of a linear fractional
transformations gives a B\"acklund or Darboux transformation.  Moreover, the orbit
of the rational loop group at the vacuum can be computed by explicit formulas. 
\item {(3)} Fordy and Kulish 
showed in 1983 ([FK]) that the MNLS equation has a sequence of commuting
Hamiltonians $\{F_j\}$ with respect to $w$ such that the Hamiltonian equation for
$F_2$ is the MNLS equation. 
\item {(4)} There is a sequence of symplectic structures
$\{w_k\}$ on \hfill\break $\cs(R,\cm_{k\times (n-k)})$:   
$$(w_k)_u(v_1, v_2)=\int_{-\infty}^\infty
\tr((W_k)_u^{-1}(v_1)v_2)dx,$$ where $W_0= -\ad(a)$ and $(W_k)_u$
is an integro-differential operator of order $k$ (cf. [Te]).   Moreover,  there is a
Lenard-Magri relation:
$$u_t= [\K F_j, a] =W_{0}(\K F_j) = W_1(\K F_{j-1})=\cdots = (W_k)_u(\K
F_{j-k}).$$ Or equivalently, the Hamiltonian equation for $F_j$ with respect
to $w_{0}=w$ is the Hamiltonian equation for $F_{j-k}$ with respect
to $w_k$.  
\item {(5)}   Let  $H_+$ denote the group generated by 
holomorphic maps
$e_{b,j}(\l)=e^{b\l^j}$ with $b\in u(n)_a$ and $j$ a positive integer. 
We proved in [TU1] that there is an action of $H_+$ on
$(\cs(R,u(n)_a^\perp),w)$ such that:
\itemitem {(i)} the $H_+$-action is Poisson.
\itemitem {(ii)} Let $\xi_{b,j}$ denote the vector field generated by the
action of the one-parameter subgroup $e_{bt,j}$ of $H_+$.  Then
the flow  generated by $\xi_{a,2}$ is the MNLS flow,  and
$[\xi_{b_1,j_1}, \xi_{b_2,j_2}]=\xi_{[b_1,b_2],j_1+j_2}$.  In particular,
$\xi_{b,j}$ commutes with $\xi_{a,2}$.  
\itemitem {(iii)} The
flows generated by the one parameter subgroups $e_{at,j}$ are the
commuting Hamiltonian flows found by Fordy and Kulish.  

\ms

The sequence of symplectic structures $w_k$ on
$\cs(R,u(n)_a^\perp)$ can be constructed using coadjoint orbits
 of certain loop algebra of an affine algebra.
  We will use a similar method  to construct  a sequence of symplectic
structures $\tau_k$ on $C(R,Gr(k,n))$ such that 
$\tau_0=\hat \tau$ (defined by \refbz{}),
$\tau_k$ is of order  $k$ and the GNLS is Hamiltonian with respect
to $\tau_k$ for all $k$.  Here $Gr(k,n)$ is isometrically embedded as the Adjoint
orbit in $u(n)$ at $a$. 

The gauge equivalence of  Lax pairs of GNLS and MNLS equation gives rise to an
isomorphism $\Phi$ from $C(R,Gr(k,n))$ to $\cs(R, u(n)_a^\perp)$
such that the GNLS flow corresponds to the MNLS flow.  If the
isomorphism $\Phi$ were a symplectic morphism from
 $(C(R, Gr(k,n)), \hat \tau)$ to $(\cs(R, \cm_{k\times (n-k)}), w)$,
then we could translate the known Hamiltonian theory for MNLS to GNLS. 
But  $\Phi$ is not symplectic.  However, we show 
that $\Phi^*(w)$ is an order $-2$  symplectic form on
$C(R,Gr(k,n))$. In general  the pull back of the order
$k$ symplectic form $w_k$ is the order
$k-2$ symplectic form $\tau_{k-2}$ under $\Phi$.    Therefore the above properties
(1)--(5) hold for the GNLS with some minor changes.  

Every compact Hermitian symmetric space $G/K$ has a standard 
isometric embedding as an adjoint orbit in $\cg$.   So the above results should
hold for any compact Hermitian symmetric space.  The non-compact examples may
give different interesting phenomena. 

This paper is organized as follows. In section 2, we use a standard embedding
of $Gr(k,n)$ in $u(n)$ to write down the GNLS and its Lax pair, and show that
its Lax pair is gauge equivalent to that of the MNLS equation.  In section 3, we
construct the isomorphism $\Phi$ from the phase space of the GNLS to
the phase space of MNLS and compute the differential of
$\Phi$. In section 4, we show that the commuting flows of the MNLS equation
give rise to a sequence of commuting, geometric, Hamiltonian flows on the path
space of $Gr(k,n)$ with respect to $\hat\tau$.  In section 5, we use the loop
algebra of certain affine algebra to construct a sequence of symplectic structures
$\tau_k$ (of order $k$) on the phase space of the GNLS and 
 show that they are the pull back of the known  order $k+2$ symplectic
form on the phase space of the MNLS equation under $\Phi$.  We also give
a Lenard-Magri relation for the GNLS.  We use the same method to construct the
two standard symplectic forms for the KdV equation in section 6.  

Some parts of the research in this paper were carried out while the first author
was  a member of
IAS in 1997-98, and while the second author was the Distinguished Visiting
Professor of IAS in 1997-98.   We would like to thank IAS for its
generous support.

\bs

\newsection Schr\"odinger flow on $Gr(k,n)$.\par

We use a standard embedding of $Gr(k,n)$ to write down the geometric
non-linear Schr\"odinger equation (GNLS) and its Lax pair.  We show that
the Lax pairs of GNLS and MNLS are gauge equivalent.

\refclaim[ak] Proposition. Let $<x,y>=-\tr(xy)$ denote the invariant
inner product on $u(n)$, and $M$ the Adjoint $U(n)$-orbit at
$a=\pmatrix{{i\over 2}I_k &0\cr 0& -{i\over 2}I_{n-k}\cr}$ in $u(n)$.
Then  $M$ equipped with the induced metric is the Hermitian
symmetric space Gr$(k,C^n)$, and the GNLS  with target $Gr(k,n)$ is
\refeq[af]$$\g_t=[\g, \g_{xx}], \qquad \g(x,t)\in M.$$

\proof Since the centralizer $U(n)_a=U(k)\times U(n-k)$, $M$ is diffeomorphic
to
$U(n)/(U(k)\times U(n-k))\simeq Gr(k,C^n)$. Note that both the induced
metric and the standard K\"ahler metric on $M$ are invariant under $U(n)$.  A
direct computation shows that they agree at $a$.  Hence
$M$ is the Hermitian symmetric space $Gr(k,n)$.  The complex structure of
$M$ at $y$ is given by $\ad(y)$.  Hence the symplectic form $\hat\tau$ on
$C(R,M)$ is
$$\hat\tau_\g(\d_1\g, \d_2\g)=\int_{-\infty}^\infty
<\ad(\g)(\d_1\g), \d_2\g>dx.$$

The gradient
of the energy functional
$E$ is $-\K_{\g_x}\g_x$, where $\K$ is the covariant derivative.  So the
Hamiltonian equation for $E$ (the GNLS on $Gr(k,n)$) is 
$$\g_t= [ -\K_{\g_x}\g_x,\g]=[\g, \K_{\g_x}\g_x].$$
Let $TM$ and $\nu(M)$ denote the tangent and normal bundles of $M$ in $u(n)$,
and $\pi_t$, $\pi_n$ the orthogonal projection onto $TM$ and $\nu(M)$
respectively.  
 Since the metric on $M$ is the induced metric, 
$$\K_{\g_x}\g_x = \pi_t(\g_{xx})= r_{xx}-\pi_n(\g_{xx}).$$
But if $\g=gag^{-1}$, then $\nu(M)_\g=gu(n)_a g^{-1}$.  So 
$[\g, v]=0$ for all $v\in \nu(M)_\g$.  This implies that
$[\g,\pi_n(\g_{xx})]=0$.  
 Hence the GNLS is equation \refaf{}. \qed

Next, we give
a Lax pair for equation \refaf{}. 

\refclaim[al] Proposition. Let $M$ be the adjoint orbit at $a$ as in
Proposition \refak{}. Then $\g$ satisfies equation \refaf{} if and only if 
\refeq[ag]$$\left[{\p\over \p x}+ \g \l, \quad {\p\over \p t} + \g \l^2 +
[\g,\g_x]\l\right]=0  \quad {\rm for\,\, all\,}\,\,\l\in C.$$

\proof  Equation \refag{} is equivalent to 
\refeq[ai]$$(\g\l)_t-(\g\l^2 + [\g,\g_x]\l)_x = \left[\g\l, \,\g\l^2 +
[\g,\g_x]\l\right].$$
Compare coefficient of $\l^j$ in equation \refai{} to give
$$\eqalign{&\g_t-[\g,\g_x]_x=0, \cr & \g_x +
\ad(\g)^2(\g_x)=0.\cr}$$ The first equation gives \refaf{}. The
second equation is always true because $\ad(\g)^2=-\id$ on $TM_\g$ for $\g\in
M$.  \qed

 The MNLS \refbj{} has a Lax pair
\refeq[ah]$$\left[{\p\over \p x} + a\l + u, \quad {\p\over \p t} + a\l^2 +
u\l + Q_2\right]=0,$$ where $u, Q_2$ are given by \refdv{}.

\refclaim[bo] Proposition.  The Lax pairs \refag{} and  \refah{} are
gauge equivalent. In particular, if $u$ is a solution of the MNLS \refbj{} then
$\g=EaE^{-1}$ is a solution of the GNLS, where $E$ is a solution of
$$\cases{E^{-1}E_x=u,&\cr E^{-1}E_t= Q_2.&\cr}$$

\proof Let $u$ be a solution of the MNLS equation.  Substitute $\l=0$ in the
corresponding Lax pair \refah{} to get
$\left[{\p\over \p x} + u, {\p\over \p t} + Q_2\right]=0$.
So there exists  $E$ satisfies
$$E^{-1}E_x=u, \quad E^{-1}E_t= Q_2, \quad E(0,0)=I.$$ Set
$\g=EaE^{-1}$.  Apply gauge transformation of $E$ to the Lax pair
\refah{} to get 
$$\eqalign{&E\left({\p\over \p x} + a\l + u\right)E^{-1} =
{\p\over\p x}+\g\l + EuE^{-1} - E_xE^{-1}= {\p\over \p x} + \g
\l,\cr &E\left({\p\over \p t}+ a\l^2+ u\l + Q_2\right)E^{-1} 
 = {\p\over \p t} + \g \l^2 + EuE^{-1}\l + EQ_2 E^{-1} - E_t E^{-1},\cr
&\qquad={\p\over \p t} + \g \l^2 + EuE^{-1}\l.
\cr}$$ 
Since $\g_x= E[u,a]E^{-1}$, $[\g,\g_x]= E[a,
[u,a]]E^{-1}$.  But
$\ad(a)^2$ is $-$id on $u(n)_a^\perp$. Hence $[\g,\g_x]=EuE^{-1}$.   This
proves that the gauge transformation of  the Lax pair
\refah{} by $E$ is the Lax pair \refag{}.  \qed

\ms

\newsection Development map of  an Adjoint orbit.\par

Let $a\in u(n)$ be a fixed diagonal matrix, and $M$ the adjoint
$U(n)$-orbit at $a$ in $u(n)$. So $M$ is a flag
manifold if $a$ has distinct eigenvalues, and is a partial flag
manifold if $a$ has multiple roots.   Let
$C_a(R,M)$ denote the space of smooth paths
$\g:R\to M$ such that 
$\g(-\infty)=\lim_{x\to -\infty} \g(x)= a$ and $\g_x:R\to u(n)$ is in the
Schwartz class.   We are motivated by the gauge equivalence given in Proposition
\refbo{} to construct an isomorphism
$\Phi$ from $C_a(R,M)$ onto the linear space  $\cs(R,u(n)_a^\perp)$.  We show
that the GNLS and MNLS flows correspond under $\Phi$.  We also will use this
isomorphism to study the relation between the hierarchies of symplectic structures
for the GNLS and MNLS flows in section 5.

\refclaim[ab] Theorem.   Let $M$ be the Adjoint orbit at $a\in u(n)$,
and  $$\Psi:\cs(R,u(n)_a^\perp)\to C_a(R,M) $$  the
map defined by $\Psi(u)= gag^{-1}$, where $g$ is the solution of
$g^{-1}g_x= u$ such that 
$g(-\infty)= I$. Then $\Psi$ is an isomorphism. 

\proof  First we prove that $\Psi$ is one to one. Assume
$\Psi(f)=\Psi(g)$. Then we have
$$\eqalign{& faf^{-1}=gag^{-1}, \quad f^{-1}f_x\in
u(n)_a^\perp,\quad g^{-1}g_x\in u(n)_a^\perp\cr &f(-\infty)=
g(-\infty)= I.\cr}$$ Set $h=f^{-1}g$.  
Since $hah^{-1}=a$,  $h(x)\in U(n)_a$ for all $x$.  Note 
$$g^{-1}g_x = (fh)^{-1}(fh)_x = h^{-1} f^{-1}f_x h + h^{-1}h_x.$$  But
$h^{-1} u(n)_a^\perp\, h\, \subset\, u(n)_a^\perp$ and
$g^{-1}g_x\in u(n)_a^\perp$ imply that
$h^{-1}h_x$ lies in $u(n)_a^\perp$. But $h^{-1}h_x$ also lies in
$u(n)_a$.  So $h^{-1}h_x= 0$, which  implies that $h$ is
constant.  But $h(-\infty) =I$.  So
$h=I$ and $f=g$. This proves $\Psi$ is one to one. 

To prove $\Psi$ is onto, given $\g\in C_a(R,M)$, choose a
smooth map
$f:R\to U(n)$ such that $\g=faf^{-1}$ and
$f(-\infty)=I$. Set $v=f^{-1}f_x$. Let $v_0$ and
$v_1$ denote the projections of $v$ onto $u(n)_a$ and
$u(n)_a^\perp$ respectively. Then $$\g_x= f[f^{-1} f_x,
a]f^{-1}=f[v,a]f^{-1} =f[v_1,a]f^{-1}.$$ Since $\g_x$  is in the Schwartz class and 
the image of
$f$ lies in the compact set $U(n)$,  $v_1$ is in the Schwartz
class.  Now solve $h:R\to U(n)_a$ such that
$$h_x=-v_0 h, \qquad h(-\infty)= I.$$ Set
$g=fh$. Then
$\g = faf^{-1} = fhah^{-1}f^{-1} = gag^{-1}$. But  
$$\eqalign{g^{-1}g_x &=h^{-1}f^{-1}f_x h + h^{-1}h_x= h^{-1}vh+ h^{-1} h_x \cr
&= h^{-1}vh - h^{-1}v_0h =h^{-1}(v-v_0)h\cr &= h^{-1}
v_1 h\quad \in\,\, h^{-1}u(n)_a^\perp h=u(n)_a^\perp . \cr}$$ So
$\Psi(g)=\g$ and
$\Psi$ is onto.

It is not difficult to see that in the appropriate topologies, the map $\Psi$ is
smooth. We explicitly compute the differential in Proposition 3.7.
\qed

\ms

\refpar[cj] Definition. The {\it development map\/} $\Phi$ on the Adjoint orbit
$M$ of $u(n)$ at $a$ is defined to be the inverse of $\Psi$, i.e.,
 $\Phi=\Psi^{-1} : C_a(R,M)\to \cs(R,u(n)_a^\perp)$. 

\ms

As a consequence of the proof of Theorem \refab{}, we have a description of
$\Phi$:

\refclaim[dm] Corollary.  If $\g\in C_a(R,M)$, then there is a unique $g:R\to U(n)$
such that $\g=gag^{-1}$, $g(-\infty)=I$, and $g^{-1}g_x\in u(n)_a^\perp$. 
Moreover, $\Phi(\g)= g^{-1}g_x$. 

We need the following Poisson operator defined in [Te] to write down the
differential of $\Phi$. A Poisson operator is a map on the tangent bundle of a
manifold, which can be used to define a Poisson or symplectic structure.   

\refpar[av] Definition.  Given $u\in \cs(R,u(n)_a^\perp)$, we introduce the
integro-differential operator $P_u:\cs(R,u(n)_a^\perp)\to C^\infty(R,
u(n)_a^\perp)$ defined by 
$$P_u(v) = v_x +\pi_a^\perp([u,v]) -[u, T_u(v)],$$ 
where $T_u(v)(x)=\int_{-\infty}^x
\pi_a([u(y),v(y)])dy$ and $\pi_a$, $\pi_a^\perp$ are the orthogonal projections of
$u(n)$ onto
$u(n)_a$ and $u(n)_a^\perp$ respectively. 
  
\ms

Note that 
$$P_u(v)=(v-T_u(v))_x + [u, v-T_u(v)]$$ and $P_u(v)$ lies
in $u(n)_a^\perp$.  So  we get a characterization of $P_u$:

\refclaim[aw] Proposition ([Te]). Let $u, v\in \cs(R,u(n)_a^\perp)$.
Then there exists a unique 
$\tilde v:R\to u(n)$ satisfies the following conditions:
\item{(i)} $\pi_a^\perp(\tilde v)=v$, where $\pi_a^\perp$ is the projection onto
$u(n)_a^\perp$,
\item {(ii)} $(\tilde v)_x + [u, \tilde v] \in u(n)_a^\perp$, 
\item {(iii)} $\tilde v(-\infty) = 0$.\ei
 Moreover,  $P_u(v) =
[d_x+ u,\tilde v]$.  

\refclaim[ci] Corollary ([Te]).  The operator $P_u$ is injective.  

\proof Suppose $P_u(v)=0$.  Let $\tilde v=v-T_u(v)$, where $T_u$ is
the operator given in Definition \refav{}.  So $P_u(v)= [d_x+u,\tilde
v]=0$, which implies that
$\tilde v(x)$ is conjugate to
$\tilde v(0)$ for all $x\in R$.  Since the conjugate class of $v(0)$ in
$u(n)$ is compact and  $\tilde v(-\infty)=0$, $\tilde v=0$. 
\qed
\ms

Recall that $\ad(a)$ is a linear isomorphism on $u(n)_a^\perp$. 

\refclaim[cx] Proposition.  Let $\g\in C_a(R,M)$, and $g:R\to U(n)$ such that
$\g=gag^{-1}$, $g^{-1}g_x= \Phi(\g)=u$ and $g(-\infty)=I$.  If
$v\in \cs(R,u(n)_a^\perp)$, then 
$$d\Phi_\g(gvg^{-1})= -P_u\ad(a)^{-1}(v).$$ 

\proof
Take variations of $\g=gag^{-1}$ and $g^{-1}g_x = u$ to get
$$ \d \g= g[g^{-1}\d g, a] g^{-1},\quad \d u= [d_x + u, g^{-1}\d g].$$ 
Set $v=[g^{-1}\d g, a]$. It follows from Proposition \refaw{} and the formula for
$\d u$ that  
$$d\Phi_\g(\d \g) = \d u = -P_u(\ad(a)^{-1} (v)). \,\qed$$

The following operators are needed later for the constructions of the
symplectic structures on $C_a(R,M)$.  

\refpar[bv]Definition. Given $\g\in C_a(R,M)$, let 
$L_\g:TC_a(R,M)_\g\to TC_a(R,M)_\g$  denote the operator defined as
follows: Write
$\g=gag^{-1}$ such that
$g^{-1}g_x=\Phi(\g)=u$ and $g(-\infty)=I$.  Then
$$L_\g(gvg^{-1})= gP_u(v)g^{-1}.$$

We give a geometric description of this operator:

\refclaim[aq] Proposition.  Let $M$ be the Adjoint orbit at $a$ in $u(n)$, and
$\eta$ a unique vector field along $\g\in C_a(R,M)$ that is tangent to $M$.  Then
there exists a unique vector field $\zeta$ along $\g$ normal to $M$ such that 
$(\eta+\zeta)_x$ is tangent to $M$ and $\zeta(-\infty)=0$. Moreover, 
$L_\g(\eta)= (\eta+\zeta)_x$.

\refclaim[bp] Theorem. Let $a={i\over 2}\pmatrix{I_k&0\cr 0& -I_{n-k}\cr}$, and
$M$ the adjoint orbit at $a$. Then the GNLS flow on
$C_a(R,M)$ corresponds to the MNLS flow on $\cs(R, u(n)_a^\perp)$
under the  development map $\Phi$.

\proof Suppose $\g(x,t)$ is a solution of the GNLS.  
 Write $\g=gag^{-1}$ such that $g^{-1}g_x\in u(n)_a^\perp$ and $\lim_{x\to
-\infty}g(x,t)=I$ for all $t$. Set $u(x,t)= g(x,t)ag(x,t)^{-1}$, i.e., $$u(\cdot, t)=
\Phi(\g(\cdot, t)).$$ Then
\refeq[dn]$$u_t= d\Phi_\g(\g_t).$$ A direct computation gives
$$\eqalign{\g_x&= g[g^{-1}g_x, a]g^{-1}= g[u,a]g^{-1},\cr
\g_{xx}&= g([u_x,a] + [u, [u,a]])g^{-1}.\cr}$$ 
Recall that 
$$\eqalign{\ck=u(n)_a&=\{y\in u(n)\n [y,a]=0\}\cr &=\left\{\pmatrix{A_1&0\cr
0&A_2\cr}\biggl| \,\, A_1\in u(k), A_2\in u(n-k)\right\},
\cr  \cp=u(n)_a^\perp &=\{y\in u(n)\n <y,z>=0 \,\, {\rm for \, all\,}\,\, z\in
u(n)_a\}\cr &=\left\{\pmatrix{0&q\cr -q^*&0\cr}\bigg| q\in \cm_{k\times
(n-k)}\right\}.\cr}$$
It follows from a simple computation that 
$$[\ck,\cp]\subset \cp, \quad [\ck, \ck]\subset \ck, \quad [\cp, \cp]\subset \ck.$$
(This is true for a Cartan decomposition of any symmetric space).  Since $[a,\ck]=0$ and
$[u,[u,a]]\in [\cp,[\cp,\ck]]\subset \ck$, we  have $$[\g,\g_{xx}]= g[a,
[u_x,a]+[u,[u,a]]]g^{-1}= g[a,[u_x,a]]g^{-1}.$$ But $\ad(a)^2= -$ id on $\cp$. So
we get
$[\g,\g_{xx}]= gu_xg^{-1}$. Substitute this into \refdn{} to get
$$u_t=d\Phi_\g(\g_t) 
=d\Phi_\g([\g, \g_{xx}])= d\Phi_\g(gu_xg^{-1}).$$
By Proposition \refcx{}, we conclude
$$u_t= -P_u\ad(a)^{-1}(u_x).$$ A direct computation shows that the right hand side is
equal to $(Q_2)_x+[u, Q_2]$, i.e., $u$ is a solution of the MNLS \refbj{}. \qed

\bs

\newsection  Commuting Hamiltonians for GNLS.\par 

We recall the definition of the orbit symplectic structure on an adjoint orbit of a
compact Lie algebra $\cu$, which is identified with $\cu^*$ via the Killing form
(trace).  Given $x\in \cu$, let $\tilde x$
denote the infinitesimal vector field generated by $x$ under the coadjoint
action.  Then the Kostant-Kirillov orbit symplectic form $\tau$ on a coadjoint
orbit $\co$ of $\cu^*$ is given by
$$\tau_\ell(\tilde x(\ell), \tilde y(\ell))= \ell([x,y]), \quad 
\ell\in \co, x, y\in \cu.$$ 
Identify $x\in \cu$ with $\ell_x\in \cu^*$ defined by $\ell_x(y)=\tr(xy)$.  Then adjoint
orbits are identified as coadjoint orbits. Let $a\in u(n)$ be a diagonal matrix, and
$M_a$ the adjoint orbit at
$a$.  Then the orbit symplectic form on $M_a$ is 
$$\tau_x(v_1, v_2)= <-\ad(x)^{-1}(v_1), v_2>.$$
It  induces a symplectic form $\hat \tau$ on the space $C_a(R,M_a)$ of maps to
$M_a$:
$$\hat\tau_\g(v_1, v_2)=\int_{-\infty}^\infty <-\ad(\g)^{-1}(v_1), v_2> dx.$$
 
Recall that $\Phi$ is the development map defined in
section 3 and $w$ is the symplectic form on $\cs(R,u(n)_a^\perp)$ defined by 
\refeq[bh]$$w(v_1, v_2)=\int_{-\infty}^\infty<-\ad(a)^{-1}(v_1), v_2>dx.$$
Then a direct computation we do later shows that $\Phi^*(w)$ is not equal to
$\tau$. We now have two symplectic
structures on $C_a(R,M_a)$. 
The flows on $C_a(R,M_a)$ pulled back from flows in the
$u(n)$-hierarchy are of course Hamiltonian with respect to $\Phi^*(w)$.
But we will show that these flows are also Hamiltonian with respect to $\hat \tau$.  

Next we recall the construction of the hierarchy for the MNLS.
Let $Q_j(u)$ be the  sequence of
 operators in $u:R\to u(n)_a^\perp$ determined by the following recursive formula
\refeq[be]$$\eqalign{&(Q_j(u))_x + [u, Q_j(u)] = [Q_{j+1}(u),a],\cr
&Q_0(u)=a,  \quad  Q_j(u)(-\infty)=0.\cr}$$
  Sattinger proved
in [Sa] (cf. also [TU1]) that $Q_j(u)$ is an order $j-1$ polynomial differential
operator in $u$. The $j$-th flow in the $u(n)$-hierarchy on
$\cs(R,u(n)_a^\perp)$ is 
\refeq[bd]$$u_t=(Q_j(u))_x + [u, Q_j(u)]= [Q_{j+1}(u),a],$$ 
which is the Hamiltonian flow for   
\refeq[dw]$$F_j(u)=-{1\over j+1} \int_{-\infty}^\infty\tr(Q_{j+2}(u)a)dx$$
with respect to $w$. Or in other words, 
\refeq[cf]$$\K F_j (u)= \pi_a^\perp(Q_{j+1}(u)),$$ 
where $\pi_a, \pi_a^\perp$ denote the orthogonal projections of $u(n)$ onto
$u(n)_a$ and $u(n)_a^\perp$ respectively.  
It follows from Proposition \refaw{} and formula \refbe{} that
\refeq[du]$$P_u(\pi_a^\perp(Q_j))= [Q_{j+1},a].$$
Or equivalently, 
\refeq[db]$$\pi_a^\perp(Q_{j+1})= -\ad(a)^{-1}P_u(\pi_a^\perp(Q_j)).$$ 
Set
\refeq[cg]$$H_j(\g)=F_j(\Phi(\g))=-{1\over j+1}\int_{-\infty}^\infty
\tr(Q_{j+2}(\Phi(\g))a)dx.$$

\refclaim[cn] Proposition. Let $\g=gag^{-1}$, $g(-\infty)=I$, and
$g^{-1}g_x=u=\Phi(\g)$. Then  $$\K H_j(\g)=
g\pi_a^\perp(Q_{j+2}(u))g^{-1}.$$

\proof Write $\d\g=gvg^{-1}$.  Then
$$\eqalign{dH_\g(\d\g)& = <\K H_j(\g), \d\g>=<g^{-1}\K H_j(\g)g, v>\cr
&=d(F_j)_u(d\Phi_\g(\d\g))= <\K F_j(u), d\Phi_\g(\d\g)>, \,\, {\rm by\,\,
\refcf{}}\cr &= <\pi_a^\perp(Q_{j+1}(u)), d\Phi_\g(\d\g)>.\cr}$$
 Proposition \refcx{} implies that
$$d\Phi_\g(gvg^{-1})= -P_u(\ad(a)^{-1}(v)).$$  Continue the above calculation,
we get
$$\eqalign{&< \pi_a^\perp(Q_{j+1}(u)), d\Phi_\g(\d\g)> =
<\pi_a^\perp(Q_{j+1}(u)), -P_u\ad(a)^{-1}(v)>
\cr & =<-\ad(a)^{-1}P_u(\pi_a^\perp(Q_{j+1}(u))), v>, \quad {\rm by\,\,}
\refdb{}\cr & =
<\pi_a^\perp(Q_{j+2}(u)), v> =<g\pi_a^\perp(Q_{j+2}(u))g^{-1},
gvg^{-1}>\cr
&=<g\pi_a^\perp(Q_{j+2}(u))g^{-1}, \d \g>.\cr}$$ 
So $\K H_j(\g)=
g\pi_a^\perp(Q_{j+2}(u))g^{-1}$.
\qed

\refclaim[dc] Corollary.  The Hamiltonian equation of $H_j=F_j\circ \Phi$
($j\geq 0$) with respect to the symplectic structure $\hat\tau$ on $C_a(R,M)$
is 
\refeq[dd]$$\g_t= g[Q_{j+2}(u),
a]g^{-1},$$ where $u=\Phi(\g)$ and $g$ satisfies $g^{-1}g_x=u$ and
$g(-\infty)=I$.  

\proof A direct computation gives
$$\eqalign{\g_t &= [\K H_j(\g), \g]\cr
&= [g\pi_a^\perp(Q_{j+2}(u))g^{-1}, gag^{-1}]\cr
&= g^{-1}[Q_{j+2}(u), a]g. \, \qed\cr}$$

\refclaim[cy] Corollary.  $\Phi$ maps the Hamiltonian flow of 
$H_{j}=F_{j}\circ \Phi$ with respect to $\hat\tau$ (equation \refdd{})
to the Hamiltonian flow of $F_{j+2}$ with respect to $w$ (the $(j+2)$-th flow
\refbd{}).

\proof  Let $\g(x,t)$ be a solution of \refdd{}, and $u(\cdot, t)= \Phi(\g(\cdot,
t))$.  Write $\g(x,t)= g(x,t)ag(x,t)^{-1}$ such that $g(-\infty, t)= I$,
$g^{-1}g_x=u$.  Then
$$\eqalign{u_t &= d\Phi(\g_t)\cr
& =d\Phi(g[Q_{j+2}(u), a]g^{-1}), \quad {\rm by\, Proposition \,}\, \refcx{},\cr
&= -P_u\ad(a)^{-1}([Q_{j+2}(u), a])\cr
&= -P_u(\pi_a^\perp(Q_{j+2}(u)), \quad {\rm by\,\,} \refdu{},\cr
&=[Q_{j+3}(u), a]. \,\qed\cr}$$

\ms

Next we show that when the adjoint orbit is the Grassman manifold  both  $H_j$
and its Hamiltonian flows are {\it geometric\/}, i.e., they only depend on the intrinsic
geometry of $Gr(k,n)$.  To show this, it suffices to show that $Q_j(\Phi(\g))$ can
be expressed in terms of covariant derivatives and the complex structure.

\refclaim[cb] Proposition. If $a={i\over 2}\pmatrix{I_k&0\cr
0& -I_{n-k}\cr}$, then $Q_j(\Phi(\g))$ is  a
polynomial in  
$\{\K_{\g_x}^k\g_x, J_\g(\K_{\g_x}^k\g_x)\n 0\leq k\leq j-1\}$.

\proof 
Since $Q_j(u)$ is a polynomial
differential operator of order $j-1$, there is a polynomial map $f_j$ such that
$$Q_j(u)= f_j(u, d^1_x u, d^2_x u, \cdots, d^{j-1}_x u),$$
where $d^k_xu = d^k u/dx^k$.  
But $\tr(gBg^{-1})=\tr(B)$, and 
$$gf_j(u, d^1_x u,  \cdots, d^{j-1}_x u)g^{-1}=  f_j(gug^{-1}, gd^1_x ug^{-1}, 
\cdots, gd^{j-1}_x ug^{-1}).$$   We claim that 
\refeq[cr]$$g d_x^ku g^{-1}=J_\g \K^k_{\g_x}(\g_x).$$  To prove this, we use
again the fact that  
$u(n)=u(n)_a + u(n)_a^\perp=\ck+\cp$ is a Cartan decomposition and 
$$TM_{gag^{-1}}=g\cp g^{-1}, \quad \nu(M)_{gag^{-1}}= g\ck g^{-1}.$$  Suppose
$u=\Phi(\g)$, $g^{-1}g_x=u$ and $g(-\infty)=I$.  Then 
$\g=gag^{-1}$ and $\g_x=g[u,a]g^{-1}$. The complex structure on  $M$ at $\g$
is $J_\g=\ad(\g)$.   Since $J_\g^2=-I$, 
$gug^{-1}= J_\g(\g_x)$. 
But $(gug^{-1})_x= g(u_x+ [u,u])g^{-1}= gu_xg^{-1}\in TM_\g$ implies that
$gu_xg^{-1}= \K_{\g_x} J_\g(\g_x)$.  Since $M$ is K\"ahler, 
$\K J= J\K$.  So $gu_xg^{-1}= J_\g(\K_{\g_x}\g_x)$.  
Then the claim follows from $[\cp,\cp]\subset \ck$ and induction. \qed

\ms 

If $a={i\over
2}\pmatrix{I_k&0\cr 0&-I_{n-k}\cr}$, then  the first few $Q_j$'s can be computed
directly from the recursive formula \refbe{}:
\refeq[cz]$$\eqalign{&Q_1(u)=u=\pmatrix{0&q\cr -q^*&0\cr}, \quad
Q_2 = \pmatrix{-iqq^*& iq_x\cr iq^*_x & iq^*q\cr},\cr
&Q_3=\pmatrix{-qq_x^* + q_xq^*& -(q_{xx}+2qq^*q)\cr
q_{xx}^* + 2q^*qq^*& q_x^*q-q^*q_x\cr}.\cr}$$ 
Note that the second flow is the MNLS flow. 
The functional $H_j=F_j\circ \Phi$ is given by the following formulae:
$$\eqalign{H_0(\g)&={1\over 2}\int_{-\infty}^\infty \N \g_x\N^2 dx, \cr
H_1(\g)&={1\over 2}\int_{-\infty}^\infty<J_\g\K_{\g_x}\g_x,
\g_x>dx,\cr
H_2(\g)&=\int_{-\infty}^\infty{1\over 2}\N \g_x\N^2 - {1\over 4}
\N\K_{\g_x}\g_x\N^4dx.\cr}$$
The corresponding Hamiltonian flows on $C_a(R, Gr(k,n))$ with respect to
$\hat \tau$ are 
$$\eqalign{&\g_t=J_\g(\K_{\g_x}\g_x),\cr
&\g_t=(\K^2_{\g_x}\g_x + 2\g_x^3),\cr
&\g_t=J_\g\left(\K_{\g_x}^3\g_x-6\g_x)\right).
\cr}$$
Note that even powers of $\g_x$ are not tensorial, but the odd powers of
$\g_x$ are for Hermitian symmetric spaces. Clearly the above flows contain only
certain expressions in these odd tensor products and their covariant derivatives.

\bs

\newsection Symplectic structures for GNLS and MNLS.\par

A two form $\w$ on  $M$  is called a {\it weakly non-degenerate\/}  if the map 
$Y: TM\to T^\ast M$ defined by $Y(v_1)(v_2)=\w(v_1, v_2)$ is injective, and
$w$ is a {\it weak symplectic form\/} if $w$ is closed and weakly non-degenerate
(cf. [CM]).  A smooth function
$f$ on
$M$ defines a  Hamiltonian vector field 
$X_f$ if $df_x$ lies in the image of $Y_x$ for all $x\in M$. Then 
$$X_f(x) = (Y_x)^{-1}(df_x).$$ 
 The Poisson bracket for two such functions $f_1, f_2$ is given by
$$\{f_1, f_2\}(x) =\w(X_{f_1}, X_{f_2})=df_1(X_{f_2}).$$
Note that when $M$ is of finite dimension, a weak symplectic 
form is symplectic,   but when $M$ is of
infinite dimension, this is not necessarily the case.

   Given a manifold $(M,g)$, assume that we can find an embedding
$f$ from $M$ to a co-adjoint orbit $N$ in the dual $\cg^*$ of a Lie algebra $\cg$
 such that $f^\ast(\o)$ is weakly non-degenerate, where $\o$ is the
Kostant-Kirillov symplectic form on
$N$.  Since $f^*(\o)$ is always closed, 
$f^*(\o)$ is a weak symplectic form on $M$.   Weak symplectic structures of
soliton equations are often constructed this way.  The procedures can be
described as follows: 
\item {(1)} Identify the algebraic structure (i.e., the reality
condition) of the Lax pair of the soliton equation, i.e., find a suitable real 
loop subalgebra of the affine Kac-Moody algebra defined by the reality
condition.
\item {(2)} Find   natural embeddings of the phase space into
coadjoint orbits of the above double loop algebra.
\item {(3)} Prove the restrictions of 
orbit symplectic forms are weakly non-degenerate. 

\ms 
The main goal of this section is to construct a hierarchy of symplectic structures
for the GNLS equation and show that the pull back of the order $k$ symplectic
form for the MNLS equation is the order $k-2$ symplectic structure in the
hierarchy.  

We first give a short review of the general construction (cf. [Te], [TU1]).  Let $\cg$
be a Lie algebra equipped with a non-degenerate bilinear form $(\,,)$.  A triple $(\cg,
\cg_+,
\cg_-)$ is a Manin-triple with respect to $(\,,)$ if  $\cg_+,\cg_-$ are Lie subalgebras of
$\cg$ such that
$\cg=\cg_+\oplus \cg_-$ as vector spaces, $\cg_+^\perp=\cg_+$, and
$\cg_-^\perp=\cg_-$.  Let $C_s(R,\cg)$ denote the Lie algebra of smooth maps
$A$ from $R$ to $\cg$ such that $A(x)(\l)$ decays as $x\to\infty$ for each $\l$. 
Then $(,)$ induces an ad-invariant form on $C_s(R,\cg)$:
$$((A,B))=\int_{-\infty}^\infty (A(x), B(x))dx.$$
Let
$$\tilde C(R, \cg) = C_s(R, \cg) + \cg \, \hat c$$ denote the Lie algebra extension
of
$C(R,\cg)$ defined by
$$[A,B]_0 =[A,B] + \rho(A,B) \hat c.$$
Here $\hat c$ is the generator of the center of the extension and $\rho$ defined by
$$\rho(A,B)= \int_{-\infty}^\infty (A_x(x),B(x)) dx$$
is a cocycle on $C_s(R,\cg)$.  In order for the expression to be skew adjoint, we
must have decay at infinity.  This requirement will require later that we impose
constraints to insure that the higher order symplectic forms are well-defined. 

 Let $C(R,\cg)$ denote the smooth maps from $R$ to $\cg$ such that
$A(x)(\l)$ is bounded in $x$ for each $\l$.  Then the set $d_x+ C(R,\cg)$ can be
identified as a subset of the dual
$\tilde C(R,\cg)^*$ via
$$\eqalign{&d_x(c)= 1, \quad d_x(B)=0,\cr
& A(B)=\int_{-\infty}^\infty (A(x), B(x))dx,\quad A(c)=0,\cr}$$
where $A\in C(R,\cg)$ and $B\in C_s(R, \cg)$.  The coadjoint action $\ast$ of
$g\in C_s(R,G)$ at
$d_x+A$ is given by the gauge transformation:
$$g\ast (d_x+A)= g(d_x+A)g^{-1}= d_x + gAg^{-1}- g_xg^{-1}.$$
For more detail see [Ka] and [PS].  
 The set $d_x+ C(R,\cg)$ is invariant under
the coadjoint action. The phase space of a soliton equation often occurs as an
coadjoint 
$\tilde C(R,\cg_-)$-orbit, which is  equipped with the orbit symplectic form.   The
second symplectic form for the soliton equation
 is often obtained by finding a new embedding of the phase space in another
coadjoint orbit via an ad-invariant bilinear form
$(,)_\L$ on $\cg$ defined by some ``Casimir'' operator $\L$. Here  a ``Casimir''
operator on
$\cg$  is a linear isomorphism $\L$ of $\cg$ that satisfies the condition: 
$$\L([x,y])=[\L(x), y] = [x, \L(y)]$$ for all $x,y\in \cg$.  Usually $\L$ is 
multiplication by an element in the center of the enveloping algebra.  Then the
bilinear form
$( , )_\L$ defined by $\L$,
$$(x,y)_\L= (\L(x), y),$$ is  non-degenerate and ad-invariant. 

\ms

We apply the abstract construction above to the $u(n)$-hierarchy.  Let $\cg$ be the
Lie algebra of holomorphic maps $A$ from
$\infty>\n\l\n>{1\over
\e}$ to $GL(n,C)$ that satisfies the $u(n)$-reality condition
$$A(\bar\l)^*+A(\l)=0.$$ Let $\cg_+$ denote the subalgebra of $A\in \cg$ such that
$A$ can be extended holomorphically to $C$, and $\cg_-$ the subalgebra of $A\in \cg$
such that
$A$ can be extended holomorphically to  a neighborhood of $\infty$ in $S^2=C\cup
\{\infty\}$ and $A(\infty)=0$. Note that  
$$A(\l)=\cases{\sum_{j\geq 0} A_j\l^j, & if $A\in \cg_+$,\cr
\sum_{j<0} A_j\l^j, & if $A\in \cg_-$.\cr}$$
Let $< , >$ denote the ad-invariant bi-linear form 
$$<A,B>= \oint \tr(A(\l)B(\l))d\l$$
on $\cg$.  Or equivalently,  if $A(\l)= \sum_j A_j \l^j$ and $B=\sum_k B_k\l^k$,
then 
$$<A, B>=\sum_j\tr(A_jB_{-j-1}).$$
Then $(\cg, \cg_+, \cg_-)$ is a Manin triple with respect to $< , >$. 
Let $a$ be a fixed diagonal matrix in $u(n)$. A direct
computation shows that the coadjoint
$\tilde C(R, \cg_-)$-orbit $\cs_0$ at $d_x+ a\l$ is
$$\cs_0=\{d_x+a\l + u\n u\in \cs(R,u(n)_a^\perp)\},$$ and the orbit symplectic form
$w_0$ is the zero order symplectic form $w$ defined by \refbh{}. The formula is
$$w_0(A,B)=\int_{-\infty}^\infty <-\ad(a)^{-1}(A(x)), B(x)> dx.$$

For an integer $k$, let $\L_k:\cg\to \cg$ denote the operator defined by 
$$\L_k(A)(\l)= \l^{-k} A(\l).$$ Then $\L_k$ 
is a Casimir operator, and 
\refeq[dy]$$<A,B>_{\L_k}=<\L_k(A), B>= \sum_j\tr(A_iB_{-i+k-1}).$$

For $k< 0$, the space $d_x+A(x)$ such that 
$$A(x)(\l)= \sum_{j\geq k} \xi_j(x)\l^j$$ is identified as a subspace of the dual of $\tilde
C(R,\cg_-)$ via $<< , >>_{\L_k}$ and is invariant under the coadjoint action.  The
infinitesimal vector field generated by
$$\xi_-=\xi_{-1}(x)\l^{-1} + \cdots + \xi_{-j}(x)\l^{-j} +\cdots \quad \in
C_s(R,\cg_-)$$ is $$\xi_-(A)=\pi_{k,\infty}([\xi_-, d_x+A]).$$ 
Here $\pi_{k,\infty}$ is the projection 
$$\pi_{k, \infty}(\sum_jA_j\l^j)= \sum_{j\geq k}A_j\l^j.$$
Let 
$$S_k= \,{\rm the \, coadjoint \,} \tilde C(R,\cg_-)- {\rm orbit\, at \,} d_x+a\l.$$
The orbit
symplectic form $\o_k$ on the coadjoint orbit $\cs_k$ is 
\refeq[bx]$$\o_k(\d_1u, \d_2 u)= <<\xi_-, \eta_-(A)>>_{\L_k}= \int_{-\infty}^\infty
\tr(\xi_{k-1}(x)\d_2u(x))dx,$$
where $\d_1u=\xi_-(A)$ and $\d_2u= \eta_-(A)$. 

We claim that $\cs_k\cap\cs_0$ is a finite codimension submanifold of $\cs_0$. 
 If $\d u$ is tangent to $\cs_k\cap \cs_0$ at $d_x+ a\l + u$, then there exists 
$$\xi_-=\xi_{-1}(x)\l^{-1} + \xi_{-2}(x)\l^{-2} +\cdots \quad\in \cg_-$$
such that  $$\eqalign{&\d u = \xi_- (d_x+a\l + u) =\pi_{k,\infty}([\xi_-,
d_x+a\l+u])\cr &= [\xi_{-1},a] +
\l^{-1}\left([\xi_{-2}, a] + [\xi_{-1}, d_x+u]\right) + \cdots + \l^k\left([\xi_{k-1}, a] +
[\xi_k, d_x+u]\right).\cr}$$
Compare coefficients of $\l^j$ in the above equation to get
\refeq[do]$$\cases{[\xi_{-1}, a]=\d u,&\cr 
(\xi_i)_x + [u,\xi_i] = [\xi_{i-1}, a],& if $k\leq i\leq -1$,\cr
\xi_i(x)\in\,\, \cs(R, u(n)), & for all $k\leq i\leq -1$.\cr}$$ 
 Given $\d u\in \cs(R, u(n)_a^\perp)$,
we have
$$\pi_a^\perp(\xi_{-1})=-\ad(a)^{-1}(\d u)$$ and
$$\pi_a(\xi_{-1}(x))= -\int_{-\infty}^x \pi_a([u,
\pi_a^\perp(\xi_{-1})]dy.$$
 So $\xi_{-1}$ is in the Schwartz
class if $\int_{-\infty}^\infty [u, \ad(a)^{-1}(\d u)]dx =0$.  Continuing
inductively, we get a formula for the $u(n)_a^\perp$ part
\refeq[eh]$$\pi_a^\perp(\xi_{i-1})=-\ad(a)^{-1}((\xi_i)_x + [u, \xi_i]).$$
For the $u(n)_a$ part, the integral formula
$$\pi_a(\xi_i)=-\int_{-\infty}^x \pi_a([u, \pi_a^\perp(\xi_i)]dy$$
leaves with constraints
$$\int_{-\infty}^\infty \pi_a([u, \pi_a^\perp(\xi_i)]dx=0,$$
which are necessary for $\pi_a(\xi_i)$ in Schwartz class for $k\leq i\leq -1$.  
Hence $\cs_0\cap \cs_k$ is of finite
codimension in $\cs_0$.  

 Use the Poisson operator $P_u$ defined in
\refav{} to write \refeh{} as
\refeq[dp]$$\pi_a^\perp(\xi_{i-1})=
(-\ad(a)^{-1}P_u)(\pi_a^\perp(\xi_1)).$$ By induction, 
\refeq[dq]$$\pi_a^\perp(\xi_{k-1})=(-\ad(a)^{-1}P_u)^{-k}(-\ad(a)^{-1}
(\d u)).$$
Let $w_k$ denote the restriction of the orbit symplectic form $\o_k$ to
$\cs_0\cap \cs_k$.   Substituting  \refdq{} to  \refbx{}, we get for $k\leq 0$
$$(w_k)_u(\d_1 u, \d_2 u)= \int_{-\infty}^\infty
(-1)^{-k+1}\tr\left(\left((\ad(a)^{-1}P_u)^{-k}\ad(a)^{-1}(\d_1
u)\right)\d_2u\right) dx.$$
Since $\ad(a)^{-1}$ is an isomorphism and $P_u$ is injective (Proposition
\refaw{}), $w_k$ is a weak symplectic form on $\cs_0\cap \cs_k$. 

For $k=1$, the space of $d_x+ A(x)$ such that $A(x)(\l)$ is of the form 
$\sum_{j\leq
0}A_j(x)\l^j$ can be identified as a subset of $\tilde C(R, \cg_+)^*$ via
$<,>_{\L_1}$.  The coadjoint orbit $\cs_1$ at $d_x$ is a submanifold of  $d_x+
 \cs(R, u(n))$, and $\d u$ is tangent to $\cs_1$ at
$d_x+u$ if and only if there exists $\xi\in \cs(R, u(n))$ such that 
$$\d u=[d_x+u, \xi].$$ For example, $v$ is tangent to $S_1$ at $d_x$ if and only
if $v=\xi_x$ for some $\xi$ in the Schwartz class.  So $v$ must satisfies the
condition $\int_{-\infty}^\infty v(x)dx=0$. 
Let $w_1$ denote  the restriction of the orbit symplectic form of $\cs_1$ to
$\cs_0\cap \cs_1$. A direct computation shows that
$$(w_1)_u(\d_1 u, \d_2 u)=\int_{-\infty}^\infty \tr(P_u^{-1}(\d_1 u)\d_2 u) dx.$$

For $k>1$, the space of $d_x + A(x)$ such that $A(x)(\l)$ is of the
form $\sum_{j\leq k-1} A_j(x)\l^j$ can be identified as a subset of $\tilde C(R,
\cg_+)$, and it is invariant under the coadjoint action. Let $\cs_k$ denote the
coadjoint
$\tilde C(R,\cg_+)$-orbit at $d_x+a\l$, and
$\o_k$ its orbit symplectic form.   Note that $\d u$ is tangent to 
$\cs_k\cap \cs_0$ if there exist $\xi_0, \cdots, \xi_{k-1}$ such that 
$$\cases{[d_x+u, \xi_0]= -\d u,& \cr
[d_x+u, \xi_i]= [\xi_{i-1}, a], & if $0\leq i\leq k-1$, \cr
\lim_{x\to \pm \infty}\xi_i(x)=0, & if $0\leq i\leq k-1$.\cr}$$
The ordinary differential equations can always be solved.  Finite constraints appear
from requiring $\xi_i(\infty)=0$.  Let $w_k$ denote the restriction of $\o_k$ to
$\cs_0\cap
\cs_k$. The formula for $w_k$ can be computed similarly as in the case $k\leq 0$.
To summarize, we have

\refclaim[dr] Theorem.  
\item {(i)} $\cs_k\cap \cs_0$ is a finite dimensional
submanifold of $\cs_0$. 
\item {(ii)} $w_k$ is an order $k$ weak symplectic form on 
$\cs_k\cap \cs_0$, and 
$$(w_k)_u(\d_1 u, \d_2 u)= \int_{-\infty}^\infty
(-1)^{-k+1}\tr(((\ad(a)^{-1}P_u)^{-k}\ad(a)^{-1}(\d_1
u))\d_2u) dx.$$

To construct symplectic structures for the GNLS flow, we proceed in the same way
as for the MNLS flow, except that we use a different Manin triple.  Let $\cg$ be
as above, and 
$$(A,B)=\sum_i \tr(A_iB_{-i+1}).$$ Let $\cg_{>0}$ denote the subalgebra of
$A\in \cg$ that can be holomorphically extended to $C$ and $A(0)=0$, and
$\cg_{\leq 0}$ the subalgebra of $A\in \cg$ that can be holomorphically extended
to a neighborhood of
$\infty$ in
$S^2=C\cup \{\infty\}$.  Then:
\item {(i)} $(\cg, \cg_{>0}, \cg_{\leq 0})$ is a Manin triple with
respect to $(\, , )$.
\item {(ii)} $(A,B)_{\L_k}= \sum_i\tr(A_i B_{k+1-i})$.

\ms

The space $C(R, \cg_{>0})$ can be identified as a subspace of $\tilde
C(R, \cg_{\leq 0})^*$ via $((\, , ))$ and is invariant under the coadjoint action. Let
$\cm_0$ denote the coadjoint 
$\tilde C(R, \cg_{\leq 0})$-orbit at $a\l$. Then $\cm_0$ is the set of $d_x+\g\l$,
where $\g(x)$ lies in the adjoint orbit of $a$ in $U(n)$, and the orbit symplectic
form $\tau_0$ is $\hat \tau$, i.e., 
$$\hat\tau_\g(\d_1 \g, \d_2 \g)= \int_{-\infty}^\infty \tr(-\ad(\g)^{-1}(\d_1\g),
\d_2 \g) dx.$$

For $k<0$, the space of $d_x+A(x)$ so that $$A(x)(\l)=\sum_{j\geq k+1}
A_j(x)\l^j$$ can be identified as a subset of $\tilde C(R, \cg_{\leq 0})^*$ via the bilinear
form $(\, , )_{\L_k}$ and is invariant under the coadjoint action.  Let $\cm_k$
denote the coadjoint $\tilde C(R,\cg_{\leq 0})$-orbit at $d_x+a\l$, $\o_k$
the orbit symplectic form, and $\tau_k$ the restriction of $\o_k$ to
$\cm_0\cap\cm_k$. 
\ms

For $k>0$, the space of $d_x+A(x)$ so that
$$A(x)(\l)=\sum_{j\leq k} A_j(x)\l^j$$ can be identified as a subset of $\tilde C(R,
\cg_{> 0})^*$ via the bilinear form $(\, , )_{\L_k}$ and is invariant under the coadjoint
action. Let $\cm_k$ denote
the coadjoint $\tilde C(R,\cg_{> 0})$-orbit at $d_x+a\l$, $\o_k$ the orbit
symplectic form, and $\tau_k$ the restriction of $\o_k$ to $\cm_0\cap \cm_k$. 
We go into less detail here, as the construction and proof are small modifications of
the MNLS case. 

\refclaim[ar] Theorem. (i) $\cm_0\cap \cm_k$ is a finite
codimension submanifold of $C_a(R, M_a)$,
\item {(ii)} $\tau_k$ is an order $k$ weak symplectic form on $\cm_0\cap\cm_k$
and 
\refeq[ax]$$(\tau_k)_\g(\d_1\g,
\d_2\g)=\int_{-\infty}^\infty
(-1)^{-k+1}\tr(((\ad(\g)^{-1}L_\g)^{-k}\ad(\g)^{-1}(\d_1\g))
\d_2\g) dx,$$ where
$L_\g$ is the operator defined in Definition \refbv{}. 

\proof We will give the computation for $k< 0$. The computation for $k>0$ is
similar.  Note that 
$\d \g$ lies in the tangent space of
$\cm_0\cap \cm_k$ if and only if there exist $\xi_-$ so that
\refeq[am]$$(\d \g) \l = \pi_{k+1,\infty}\left(\left[\xi_-,
d_x+\g\l\right]\right),$$ where $\pi_{k+1,\infty}(\sum_j A_j\l^j)= \sum_{j\geq
k+1} A_j\l^j$.  The calculation below is entirely algebraic.   Write 
$$\xi_-=\xi_0 +\xi_{-1}\l^{-1} + \xi_{-2}\l^{-2} + \cdots.$$ 
Equation \refam{} gives
\refeq[ap]$$\cases{[\xi_0,\g]=\d \g, &\cr
(\xi_j)_x -[\xi_{j-1},\g]=0, & if $ k+1\leq j\leq 0$, \cr
\lim_{x\to \pm\infty}\xi_j(x)= 0, & if $  k+1\leq j\leq 0$.\cr}$$ 
Let $\xi_j^T$ denote the projection of the vector field $\xi_j$ along
$\g$ to $TM_\g$.  Proposition
\refaq{} implies that if $\xi_j$'s satisfy \refap{} then
\refeq[ds]$$\xi_j^T(\d \g) = (-\ad(\g)^{-1}L_\g)^{-j}(-\ad(\g)^{-1})(\d\g).$$

 Given two tangent vectors  $\d_1\g, \d_2\g$,  
there exist $\xi, \eta\in C(R,_{-\infty,0})$ such that 
$$\l\d_1\g=\pi_{k+1, \infty}([\xi,d_x+\g\l]), \quad \l \d_2 \g=
\pi_{k+1, \infty}([\eta,d_x+\g\l]).$$
Write $\xi=\sum_{j\leq 0} \xi_j \l^j$ and $\eta=\sum_{j\leq 0} \eta_j\l^j$.
Then the restriction of the orbit symplectic form on $\cm_k$ to
$\cm_k\cap\cm_0$ is
\refeq[dt]$$\eqalign{\tau_k(\d_1\g, \d_2\g) 
&=  ((\xi, \,\,(\d_2 \g)\l))_{\L_k}\cr
&= \int_{-\infty}^\infty \tr(\xi_k, \d_2\g) dx.\cr}$$
Substitute \refds{} into \refdt{} to get formula \refax{}.  
\qed

\refclaim[ca] Theorem.   $\Phi^*(w_k)=\tau_{k-2}$.

\proof  Write $\g=gag^{-1}$ such that $g^{-1}g_x=\Phi(\g)=u$ and $g(-\infty)=I$. 
Then
$$\d_i \g = g[g^{-1}\d_i g, a]g^{-1}.$$ Set 
$$v_i=[g^{-1}\d_i g, a],\quad I_a= \ad(a), \quad Y_\g= \ad(\g).$$ 
Let $\li\,\,\ri_o$ denote the $L^2$ inner product.  We compute
$\Phi^\ast(w_k)$.
$$\eqalign{(\Phi^*(w_k))_\g (\d_1 \g, \d_2\g) &=
(w_k)_u(d\Phi_\g(\d_1\g), d\Phi_\g(\d_2\g)), \quad {\rm by\,\, Proposition\,\,}
\refcx{}\cr  &= w_k(P_uI_a^{-1}(v_1),
P_uI_a^{-1}(v_2)),\quad {\rm by\,\, Theorem\,\,} \refdr{}\cr 
&=(-1)^{-k+1}\li (I_a^{-1}P_u)^{-k}I_a^{-1}P_uI_a^{-1}(v_1),
P_uI_a^{-1}(v_2)\ri_o\cr
&=(-1)^{-k+1}\li(I_a^{-1}P_u)^{-k+1}I_a^{-1}(v_1), P_uI_a^{-1}(v_2)\ri_o\cr
&=(-1)^{-k+1}\li(I_a^{-1}P_u)^{-k+2}I_a^{-1}(v_1), v_2\ri_o\cr
 &=(-1)^{-k+1}\li g((I_a^{-1}P_u)^{-k+2}I_a^{-1} v_1)g^{-1},
gv_2g^{-1}\ri_o.\cr}$$ Note that $$gI_a(v)g^{-1}= g[a,v]g^{-1}= [gag^{-1},
gvg^{-1}]=[\g, gvg^{-1}]= Y_\g(gvg^{-1}).$$  The operator
$L_\g$ is defined in \refbv{} by $L_\g(gvg^{-1})= gP_u(v)g^{-1}$.  So we have
$$Y_\g^{-1}L_\g(gvg^{-1})= g(I_a^{-1}P_u(v))g^{-1}.$$ Hence
$$\Phi^*(w_k)_\g (\d_1 \g,
\d_2\g)=(-1)^{-k+1}\li(Y_\g^{-1}L_\g)^{-k+2}Y_\g^{-1}(\d_1\g),
\d_2\g\ri_o,$$ which is equal to $(\tau_{k-2})_\g(\d_1\g, \d_2\g)$. \qed

\ms

\refclaim[co] Corollary.  If $j\geq 1$, then the  
Hamiltonian equation for $H_j=F_j\circ \Phi$ with respect to
$\tau_0=\hat\tau$ satisfies the following Lenard-Magri relation
\refeq[cp]$$\g_t=[\K H_j(\g), \g]=L_\g(\K H_{j-1}(\g)),$$  where 
$L_\g$ is defined in \refbv{}.  In other words, the Hamiltonian equation for $H_j$
with respect to $\tau_0$ is the Hamiltonian equation for $H_{j-1}$ with respect to
$\tau_1$. 

\proof  The Hamiltonian equation for $H_j$ with respect to $\tau_0$  is 
$\g_t=[\K H_j(\g), \g]$.  By 
 Proposition \refcn{}, we have
$$\eqalign{L_\g(\K H_{j-1}(\g))&= L_\g(g\pi_a^\perp(Q_{j+1}(u))g^{-1}),\cr
&=gP_u(\pi_a^\perp(Q_{j+1}(u)))g^{-1}, \quad {\rm by\,\,} \refdu{},\cr &=
g[Q_{j+2}(u), a]g^{-1}= [gQ_{j+2}(u)g^{-1}, \g] = [\K H_j(\g), \g]. \qed\cr}$$
\bs

\newsection Symplectic structures for KdV.\par

Since we have set up the machinery for constructing symplectic structures and
applied it to two examples, we take this opportunity to show that two structures
for KdV can be obtained in the same fashion.  

The KdV equation,
\refeq[eg]$$q_t={1\over  4}(q_{xxx} - 6qq_x),$$ has a  Lax pair:
$$\left[{\p\over \p x} + a\l + u, \,\, {\p\over \p t} + a\l^3 +
u\l^2 + Q_2\l + Q_3\right]=0,$$ where 
$$\eqalign{ & a= \pmatrix{1&0\cr 0&-1\cr}, \quad u=
\pmatrix{0& q\cr 1& 0\cr}, \cr
& Q_2= \pmatrix{-{q\over 2}& -{q_x\over
2}\cr 0 & {q\over 2}\cr},\quad Q_3= \pmatrix{{q_x\over 4}& {q_{xx} -
2q^2\over 4}\cr -{q\over 2}& -{q_x\over 4}\cr}.\cr}$$
This Lax pair satisfies the  following reality condition:
\refeq[hc]$$\cases{\overline{A(\bar \l)}=A(\l),& \cr
\phi(\l)^{-1}A(\l)\phi(\l) =
\phi(-\l)^{-1} A(-\l) \phi(-\l),& where  $\phi(\l)=\pmatrix{1&\l\cr 0
&1\cr}$.\cr}$$
We call this the {\it KdV reality condition\/} (cf. [TU2]). 
 
 Let $\cg, \cg_+, \cg_-, < , >, < , >_{\L_k}$ be as in section 5, and $\cg^{kdv}$
($\cg^{kdv}_\pm$ resp.) the space of all $A\in \cg$ ($\cg_\pm$ resp.) that
satisfies the KdV-reality condition.  Recall that if $\xi= \sum_i\xi_i\l^i$ and
$\eta=\sum_j \eta_j\l^j$, then the bilinear form $< , >_{\L_k}$
defined by \refdy{} is
$$<\xi,\eta>_{\L_k}=\sum_i\tr(\xi_i\eta_{-i+k-1}).$$

  Let $e_{ij}\in sl(2)$ denote the matrix whose entries
are zero except the $ij$-th entry equals to $1$. 

\refclaim[ei] Lemma.  Let $\xi(\l)=\sum_j \xi_j \l^j$ with
$\xi_j=\pmatrix{A_j &B_j\cr C_j& -A_j\cr}\in sl(2,R)$.  Then $\xi$ satisfies the
KdV-reality condition if and only if 
\refeq[ea]$$\xi_{2j}=\pmatrix{A_{2j}&B_{2j}\cr
C_{2j}&-A_{2j}\cr}, \quad\xi_{2j+1}=\pmatrix{C_{2j}& -2A_{2j}\cr 0&-C_{2j}\cr}$$ for
all
$j$. 

\proof   $\xi$ satisfies the KdV-reality condition if and only if the coefficient of
$\l^{2j+1}$ in $\phi(\l)^{-1}\xi_-(\l) \phi(\l)$ is zero for all $j$, i.e., 
\refeq[dg]$$ \cases{A_{2j+1}-C_{2j}=0, &\cr 
B_{2j+1} + 2A_{2j} - C_{2j-1}=0, &\cr 
C_{2j+1}=0.&\cr}$$ 
So we have
\refeq[dx]$$C_{2j+1}=0, \quad A_{2j+1}=C_{2j}, \quad B_{2j+1}= - 2A_{2j},$$
which proves the Lemma. \qed

\refclaim[dz] Proposition. The restriction of the bilinear form $< , >_{\L_k}$
 to $\cg^{kdv}$ is degenerate if $k$ is even, and is
non-degenerate if $k$ is odd. 

\proof  
 By Lemma \refei{}, $be_{12}\in
\cg^{kdv}$.  If $k$ is even, then 
$$<be_{12}, \xi>_{\L_k}=bC_{k-1}$$ for all  $\xi\in \cg^{kdv}$, where
$\xi=\sum_i \xi_i\l^i$ and $\xi_i=\pmatrix{A_i&B_i\cr C_i & -A_i\cr}$.
 But $k-1$ is
odd,  so
$C_{k-1}=0$. This shows that $< , >_{\L_k}$ is degenerate.  

Next we prove that if $k$ is odd then $< , >_{\L_k}$ is non-degenerate on
$\cg^{kdv}$.  Let $\xi=\sum_i\xi_i\l^i, \eta=\sum_j \eta_j\l^j\in \cg^{kdv}$,
and $$\xi_i=\pmatrix{A_i&B_i\cr C_i& -A_i\cr}, \quad
\eta_i=\pmatrix{A_i'&B_i'\cr C_i'& -A_i'\cr}.$$   Then 
$$<\xi, \eta>_{\L_k}= \sum_i\tr(\xi_i\eta_{-i+k-1}) 
=\sum_j \tr(\xi_{2j}\eta_{-2j+k-1}) + \tr(\xi_{2j+1}\eta_{-2j+k-2}).$$
Note $-2j+k-2$ is odd. By Lemma \refei{}, we get
$$<\xi, \eta>_{\L_k}
=\sum_j 2A_{2j}A'_{-2j+k-1} + C_{2j}(B'_{-2j+k-1}+2 C'_{-2j+k-3})+
B_{2j}C'_{-2j+k-1}.$$
It follows that if $<\xi,\eta>_{\L_k}=0$ for all $\eta\in \cg^{kdv}$ then $\xi=0$. This
proves $< , >_{\L_k}$ is non-degenerate. 
\qed

It follows from Lemma \refei{} that 
$$A_0= a\l+ e_{21}$$ satisfies the KdV-reality
condition. So it belongs to $\cg_+^{kdv}$.  

Let $k\leq -1$ be an odd integer. Then the set of  $d_x+A(x)$ such that 
$$A(x)(\l)= \sum_{j\geq k} A_j(x)\l^j$$ can
be identified as a subset of $\tilde C(R,\cg_-^{kdv})^*$ via $<< , >>_{\L_k}$, and
it is invariant under the coadjoint $\tilde C(R,\cg_-^{kdv})$-action. Set
$$\cases{&$\W_{k}$= the coadjoint
$\tilde C(R, \cg^{kdv}_-)$-orbit at $d_x+a\l + e_{21}$, \cr
& $\s_k$= the orbit symplectic form on $\W_k$.\cr}$$ 

The set of $d_x+ A(x)\in  d_x + C(R, \cg^{kdv})$ such that $A(x)(\l)=\sum_{i\leq
0} A_i(x)\l^i$ can be identified as a subset of
$\tilde C(R,
\cg_+^{kdv})^*$ via $<< , >>_{\L_1}$, and is invariant under the coadjoint $\tilde
C(R,\cg_+^{kdv})$-action. Set 
$$\cases{&$\W_1$= the coadjoint $\tilde C(R, \cg_+^{kdv})$-orbit at 
$d_x+ e_{21}$, \cr
&$\s_1$= the orbit symplectic form on $\W_1$.\cr}$$

By Lemma \refei{}, $d_x+a\l+ e_{21} + u$
satisfies the KdV-reality condition if and  only if 
$$u=qe_{12} = \pmatrix{0&q \cr 0&0\cr}.$$  
Now set
$$\cn_k=\cases{\{q\in \cs(R,R)\n
(d_x+A_0 + qe_{12})\in
\W_k\}& if $k\leq -1$ is odd,\cr
\{q\in \cs(R,R)\n
(d_x+e_{21} + qe_{12})\in \W_0\}& if $k=1$.\cr}$$  
 Let $\b_k$ denote the restriction of the orbit
symplectic form on $\W_k$ to
$\cn_k$ (here we identify $\cn_k$ as a subspace of $\W_k$ via
$q\mapsto d_x+A_0+ qe_{12}$ if $k\leq -1$ is odd and $q\mapsto q+ e_{21} +
qe_{12}$ if $k=1$).  
Write 
$$(\b_k)_q(\d_1 q, \d_2
q)=\int_{-\infty}^\infty ((J_k)_q^{-1}(\d_1 q))
\d_2 q dx.$$

\refclaim[di] Theorem. 
\item {(i)} $(J_{-1})_q(v)= -2 v_x$.
\item {(ii)} $(J_1)_q(v)= {1\over 2} v_{xxx} - 2 q v_x - q_x v$.

\proof  (i)
 $\d q$ lies in the tangent space of $\cn_{-1}$ at $q$ if and only if
there exists $\xi_-\in C(R, \cg^{kdv}_-)$ such that 
$$e_{12}\d q= \pi_{-1, \infty}([\xi_-, d_x+ a\l + e_{21}+ qe_{12}]),$$
where $$\pi_{-1, \infty} (\sum_j A_j\l^j)= \sum_{j\geq -1} A_j\l^j.$$
This implies that 
\refeq[df]$$\cases{[\xi_{-1},a]= e_{12}\d q,&\cr 
[d_x+ e_{21}+qe_{12}, \xi_{-1}]= [\xi_{-2},a],\cr
\lim_{x\to \pm \infty} \xi_j(x)=0, & if $j= -1, -2$.\cr}$$

Let $\d_1u=(\d_1 q)e_{12}, \d_2u= (\d_2 q) e_{12}$ be two tangent
vectors of $\cn_{-1}$ at $q$. So  there exist $\xi=\sum_{j\leq -1}
\xi_j \l^j$ and
$\eta=\sum_{j\leq -1} \eta_j\l^j$ satisfying equation
\refdf{}. Write $$\xi_j=\pmatrix{A_j&B_j\cr C_j& -A_j\cr}.$$ By definition of the
orbit symplectic form, we get
$$(\b_{-1})_q(\d_1 q, \d_2 q) = \int_{-\infty}^\infty \tr(\xi_{-2} \d_2 q
e_{12})dx =\int_{-\infty}^\infty C_{-2}\d_2 q dx.$$   
Next we solve $\d_1 q$ in terms of $C_{-2}$. The first equation in \refdf{}
implies that
\refeq[eb]$$\d_1 q= - 2B_{-1}.$$ The second equation in \refdf{} gives
\refeq[ec]$$\cases{(A_{-1})_x-B_{-1}=0,&\cr
(B_{-1})_x - 2qA_{-1} = -2B_{-2},&\cr 
2A_{-1}= 2C_{-2}.&\cr}$$
Substitute \refeb{} into \refec{} to get $\d_1 q= - 2 (C_{-2})_x$.  This proves
$J_{-1}= - 2 d_x$.  

\ss 

(ii) If $\d q$ is tangent to $\cn_1$ at
$q$, then there exist $\xi\in C(R, \cg^{kdv}_+)$ 
$$\xi(\l) =\xi_0 + \xi_1\l + \xi_2\l^2 + \cdots$$
such that 
\refeq[dh]$$\d u= (\d q) e_{12} =\pi_{-\infty, 0}([\xi, d_x+ e_{21} +
qe_{12}]) = [\xi_0, d_x+ e_{21} +
qe_{12}].$$  Write $\xi_0=\pmatrix{A&B\cr C&-A\cr}$.
Then $$(\b_1)_q(\d_1 q, \d_2 q)= \int_{-\infty}^\infty <\xi_0, \d_2 q e_{12}>dx
=\int_{-\infty}^\infty C\d_2 q dx.$$   We need to compute $\d q$ in terms  $C$. 
To do this, we equate the entries of equation \refdh{} to get
$$\cases{A_x+qC-B=0,&\cr C_x+2A=0,&\cr 2qA- B_x= \d q.&\cr}$$ 
The second equation gives $A=-C_x/2$. Substitute this to the first equation to
solve $B$ in terms of $C$. Then the last equation solves 
$$\d q= {1\over 2}C_{xxx} -q_x C - 2 q C_x,$$ which gives the formula for $J_1$. 
\qed

\bs

\Bibliography

\a //BC//Beals, R., Coifman, R.R.//Scattering and inverse scattering for
first order systems//Commun. Pure Appl. Math.//37//1984//39-90////

\a //BFPP//Burstall, F.E., Ferus, D., Pedit, F., Pinkall, U.//Harmonic tori in
symmetric spaces and commuting Hamiltonian systems on loop
algebras//Annals of Math.//138//1993//173-212//// 

\a //BG//Burstall, F.E., Guest, M.A.//Harmonic two-spheres in compact symmetric
spaces//Math. Ann.//309//1997//541-572////

\p //CSU//Chang, N., Shatah, J., Uhlenbeck, K.//Schr\"odinger
maps//////////preprint //

\b //CM//Chernoff, P., Marsden, J.//Properties of infinite dimensional Hamiltonian
system// Lecture Notes in Math., vol. 425,  Springer-Verlag, Berlin and New
York////1974//////

\a //1906 dR//Da Rios//Sul moto d'un liquido indefinito con un filetto
vorticoso di forma qualunque (on the motion of an unbounded liquid with a
vortex filament of any shape)//Rend. Circ. Mat. Palermo//22//1906//117-135////

\b //FT//Faddeev, L.D., Takhtajan, L.A.//Hamiltonian Methods in the theory
of Solitons//Springer-Verlag////1987//////

\a //FK//Fordy, A.P., Kulish, P.P.//Nonlinear Schr\"odinger equations and simple
Lie algebra//Commun. Math. Phys.//89//1983//427-443////

\a //H//Hasimoto, H.//A soliton on a vortex filament//J.
Fluid Mechanics//51//1972//477-////

\b //Ka//Kac, V.G.//Infinite Dimensional Lie Algebras//Cambridge University
Press////1985////// 

\p //LP//Langer, J., Perline, R.//Geometric realizations of
Fordy-Kulish nonlinear Schr\"odinger systems//////////preprint //

\a //M//Magri, F.//A simple model of the integrable Hamiltonian
equation//J. Math. Phys.//19//1978//1156-1162////

\b //PS//Pressley, A. and Segal, G. B.//Loop Groups//Oxford Science Publ., Clarendon Press, Oxford//
//1986//////

\a //Sa//Sattinger, D.H.//Hamiltonian hierarchies on semi-simple Lie
algebras//Stud. Appl. Math.//72//1984//65-86////

\a //Te//Terng, C.L.//Soliton equations and differential
geometry//J. Differential Geometry//45//1997//407-445////

\p //TU1//Terng, C.L., Uhlenbeck, K.//Poisson actions and scattering theory
for integrable systems//////////to appear in J. Differential Geometry survey
volume III, preprint dg-ga 9707004//

\p //TU2//Terng, C.L., Uhlenbeck, K.//B\"acklund
transformations and loop group actions////////// to appear in Comm. Pure. Appl.
Math., preprint math.DG/9805074//

\def\uhlenaddr{{\hsize=2in \vsize=1in
\hbox{\vbox{
\leftline{Karen Uhlenbeck}
\leftline{Department of Mathematics}
\leftline{The University of Texas at Austin}
\leftline{RLM8.100 Austin, Texas 78712}
\leftline{email:uhlen@math.utexas.edu}
}}}}

\def\cltaddr{{\hsize=2in \vsize=1in
\hbox{\vbox{
\leftline{Chuu-lian Terng}
\leftline{Department of Mathematics}
\leftline{Northeastern University}
\leftline{Boston, MA  02115}
\leftline{email: terng@neu.edu}
}}}}

\bigskip\bigskip
\leftline{\cltaddr \hfill \uhlenaddr}

\end